\documentclass[oneside,a4paper,12pt,reqno]{amsart}

\usepackage{amssymb}

\newtheorem{theorem}{Theorem}[section]
\newtheorem{lemma}[theorem]{Lemma}

\newtheorem{proposition}[theorem]{Proposition}

\theoremstyle{definition}

\begin{document}
\title{Mixing actions of the rationals}
 \keywords{}
\author{Richard Miles}
\email{r.miles@uea.ac.uk}
\author{Tom Ward}
\address{School of Mathematics, University of
East Anglia, Norwich NR4 7TJ, England}
\email{t.ward@uea.ac.uk}


\begin{abstract}
We study mixing properties of algebraic actions of~$\mathbb Q^d$,
showing in particular that prime mixing~$\mathbb Q^d$ actions on
connected groups are mixing of all orders, as is the case
for~$\mathbb Z^d$-actions. This is shown using a uniform result
on the solution of~$S$-unit equations in characteristic zero fields
due to Evertse, Schlickewei and Schmidt. In contrast, algebraic actions
of the
much larger group~$\mathbb Q^*$ are shown to behave quite differently,
with finite order of mixing possible on connected groups.
\end{abstract}
\maketitle

Mixing properties of~$\mathbb Z^d$-actions by automorphisms
of a compact metrizable abelian group are quite
well understood.
Roughly speaking, the
picture has three facets.
Firstly, the correspondence between such
actions and countably generated modules
over the integral group ring~$R_d=\mathbb Z[\mathbb Z^d]$
of the acting group~$\mathbb Z^d$
due to Kitchens and Schmidt~\cite{MR1036904}
allows any mixing problem to be reduced
to the case corresponding to a cyclic module
of the form~$R_d/P$ for a prime ideal~$P\subset R_d$.
Secondly, in the connected case $P\cap\mathbb Z=\{0\}$,
Schmidt and Ward~\cite{MR1193598} showed that
mixing implies mixing of all orders by relating
the mixing property to~$S$-unit
equations and exploiting a deep result of
Schlickewei on
solutions of such equations~\cite{MR1069241}
(see also~\cite{MR766298} and~\cite{MR1119694}).
Finally, in the totally disconnected case
$P\cap\mathbb Z=p\mathbb Z$ for some rational
prime~$p$, Masser~\cite{masser} has shown that the order of
mixing is determined by the mixing behaviour of
shapes, reducing the problem -- in principle --
to an algebraic one.

Our purpose here is to show how some of this
changes for algebraic
actions of infinitely generated abelian
groups. The algebra is more involved, so
for simplicity we restrict attention to the
simplest extreme examples: actions of~$\mathbb Q^{\times}_{>0}$
(isomorphic to the direct sum of countably
many copies of~$\mathbb Z$) and actions
of~$\mathbb Q^d$ (which is a torsion extension
of~$\mathbb Z^d$). These groups are
the simplest non-trivial examples chosen from the
`dual' categories of free abelian and infinitely
divisible groups in the sense of
MacLane~\cite{MR0049192}.
The algebraic difficulties mean
we cannot present the complete picture found
for~$\mathbb Z^d$-actions, and the emphasis is
partly on
revealing or suggestive examples.

Let~$\alpha$ be an action of a countable abelian
group~$\Gamma$ on a probability space~$(X,\mathcal B,\mu)$.
For a sequence~$(\gamma_n)$ in~$\Gamma$,
write~$\gamma_n\to\infty$ if for every finite set~$F\subset\Gamma$
there is an~$N$ for which~$n>N$ implies that~$\gamma_n\notin F$.
The action~$\alpha$ is said to be \emph{mixing on~$r$ sets}
if for any sets
$$
A_1,\dots,A_r\in\mathcal B,
$$
$$
\mu\left(
\alpha_{\gamma_1}A_1\cap\dots\cap\alpha_{\gamma_r}A_r
\right)
\rightarrow\mu(A_1)\cdots\mu(A_r)\mbox{ as }
\gamma_s-\gamma_t\to\infty\mbox{ for }s\neq t.
$$
The \emph{order of mixing}~$\mathcal M(\alpha)$
of~$\alpha$ is the
largest value of~$r$ for which~$\alpha$ is mixing
on~$r$ sets, and~$\alpha$ is said to be \emph{mixing
of all orders}, denoted~$\mathcal{M}(\alpha)=\infty$,
if it is mixing on~$r$ sets for all~$r$.

\section{Algebraic actions}

Just as for algebraic~$\mathbb Z^d$-actions (see
Schmidt's monograph~\cite{MR97c:28041}),
Pontryagin duality
gives a description of~$\Gamma$-actions by automorphisms
of compact abelian groups in terms of modules over the
ring~$\mathbb Z[\Gamma]$. If~$\alpha$ is a~$\Gamma$-action
by automorphisms of~$X$, then the character group~$M=\widehat{X}$
inherits the structure of a~$\mathbb Z[\Gamma]$-module
via the dual automorphisms~$\widehat{\alpha_{\gamma}}$
for~$\gamma\in\Gamma$; conversely any~$\mathbb Z[\Gamma]$-module~$M$
defines a compact abelian group~$X_M=\widehat{M}$
carrying a dual~$\Gamma$-action~$\alpha_M$.
Write~$\lambda=\lambda_X$ for the Haar measure on~$X$.

A module is called \emph{cyclic} if it is
singly generated as a module, so takes the form~$\mathbb Z[\Gamma]/I$
for some ideal~$I\subset\mathbb Z[\Gamma]$, and
the dual~$\Gamma$-action will be called
\emph{prime} (or \emph{radical})
if the module takes the form~$\mathbb Z[\Gamma]/P$
for some prime (resp. radical) ideal~$P\subset\mathbb Z[\Gamma]$.

The rings that arise here are~$R_\infty=\mathbb Z[\mathbb Q_{>0}^{\times}]$,
corresponding to actions of~$\mathbb Q_{>0}^{\times}$,
and~$R_{\mathbb Q^d}=\mathbb Z[\mathbb Q^d]$,
corresponding to actions of~$\mathbb Q^d$.
Notice that these are wildly different rings:
for example,~$R_{\infty}$ has infinite Krull dimension,
while~$R_{\mathbb Q^d}$ has Krull dimension~$d+1$.
Both are non-Noetherian rings.

\section{Actions of~$\mathbb Q^d$}

The main result of~\cite{MR1193598} says that for an
algebraic~$\mathbb Z^d$-action~$\alpha$ on a
\emph{connected} group,
$$
\mathcal{M}(\alpha)>1\implies\mathcal{M}=\infty.
$$
The same property turns out to also hold for the simplest actions
of~$\mathbb Q^d$. This is shown in Theorem~\ref{Qcase}
below, which is stated in a slightly more general setting.
The \emph{rational rank} of an abelian group is the
maximal number of elements which are linearly
independent over~$\mathbb Z$. Thus~$\mathbb Q^d$ and~$\mathbb Z^d$
have rational rank~$d$, while~$\mathbb Q^{\times}$ does not
have finite rational rank. If~$\Gamma$ has rational
rank~$d$, then~$R_{\Gamma}=\mathbb Z[\Gamma]$ has Krull dimension~$d+1$,
and may or may not be Noetherian depending on the divisibility
properties of~$\Gamma$.

\begin{theorem}\label{Qcase}
Let~$\alpha$ be an algebraic action of
a countable torsion-free group~$\Gamma$ of finite rational
rank
corresponding to a cyclic module~$R_{\Gamma}/I$
with~$I\cap\mathbb Z=\{0\}$ and~$I$ a radical ideal.
Then
$$
\mathcal{M}(\alpha)>1\implies\mathcal{M}=\infty.
$$
\end{theorem}

\begin{proof}
The result in~\cite{MR1193598} depended
on a bound for the number of
solutions to~$S$-unit equations over number
fields~\cite{MR1069241} with certain
uniformity properties; as shown in~\cite[Sect.VIII.27]{MR97c:28041}
it is enough to have a qualitative bound for
any characteristic zero field.
Here we use instead the following deep result
from~\cite{MR1923966}.

\begin{theorem}{\sc[Evertse, Schlickewei \& Schmidt]}\label{deepmagic}
Let~$\mathbb K$ be an algebraically closed field
of characteristic zero, and let~$\Gamma$ be a
finitely generated multiplicative subgroup of~$(\mathbb K^{\times})^n$
with rank~$r$. For fixed
$$
a_1,\dots,a_n\in\mathbb K^{\times},
$$
the number of solutions~$(x_1,\dots,x_n)\in\Gamma$
of the equation
$$
a_1x_1+\dots+a_nx_n=1
$$
for which no proper subsum vanishes is bounded above
by
$$
\exp\left((6n)^{3n}(r+1)\right).
$$
\end{theorem}

Returning to the proof of Theorem~\ref{Qcase},
write~$\alpha=\alpha_{R_{\Gamma}/I}$,~$X=X_{R_{\Gamma}/I}$
and~$\lambda=\lambda_X$.
Assume that the action is not mixing on~$r$ sets
for some~$r>1$. It follows that there
are measurable sets $A_1,\dots,A_r\subset X$
and a sequence of~$r$-tuples
$$
\left(
q_1^{(j)},\dots,q_r^{(j)}
\right)_{j\ge1}\in\Gamma^r
$$
such that $q_s^{(j)}-q_t^{(j)}\to\infty$ as~$j\to\infty$
for~$s\neq t$, for which
\begin{equation}\label{Qcasenotmixingonrsets}
\lambda\left(
\alpha_{q_1^{(j)}}A_1\cap\dots\cap
\alpha_{q_r^{(j)}}A_r\right)\not\rightarrow
\prod_{s=1}^{r}\lambda(A_s)\mbox{ as }j\to\infty.
\end{equation}
By approximating the indicator functions of the
sets appearing in~\eqref{Qcasenotmixingonrsets}
and applying the orthogonality relations for
characters on~$X$, it follows that there are
non-zero
elements~$a_1,\dots,a_r\in R_{\Gamma}/I$
with the property that
\begin{equation}\label{eq:Qcaseone}
\widehat{\alpha_{q_1^{(j)}}}(a_1)
+\dots+
\widehat{\alpha_{q_r^{(j)}}}(a_r)=0\mbox{ for infinitely many }j.
\end{equation}
First assume that the ideal~$I$ is a prime ideal~$P$.
Embed~$R_{\Gamma}/P$ into a field~$\mathbb K$ of
characteristic zero (this is possible because~$P$
is prime and~$P\cap\mathbb Z=\{0\}$); denote by
$$
x\mapsto
u_1^{q_{s,1}^{(j)}}\cdots
u_d^{q_{s,d}^{(j)}}x=\boldsymbol{u}^{\boldsymbol{q}_{s}^{(j)}}x
$$
the automorphism
of~$\mathbb K$ defined by the
automorphism~$\widehat{\alpha_{\boldsymbol{q}_s^{(j)}}}$ of~$R_{\Gamma}/P$
(writing~~$\boldsymbol{u}^{\boldsymbol{q}}$
for~$u_1^{q_1}\cdots u_d^{q_d}$,
where~$\boldsymbol{q}=(q_1,\dots,q_d)\in\Gamma^d$).
Then equation~\eqref{eq:Qcaseone} implies that
\begin{equation}\label{eq:Qcasetwo}
\boldsymbol{u}^{\boldsymbol{q}_1^{(j)}}a_1
+\dots+
\boldsymbol{u}^{\boldsymbol{q}_r^{(j)}}a_r=0\mbox{ for infinitely many }j
\end{equation}
holds in~$\mathbb K$.
Rearranging, this gives an equation
\begin{equation}\label{eq:Qcasethree}
(-a_2/a_1)\boldsymbol{u}^{\boldsymbol{q}_2^{(j)}-
\boldsymbol{q}_1^{(j)}}+\dots+
(-a_r/a_1)\boldsymbol{u}^{\boldsymbol{q}_r^{(j)}-
\boldsymbol{q}_1^{(j)}}=1\mbox{ for infinitely many }j.
\end{equation}
Assume that the rational rank of~$\Gamma$ is~$d<\infty$
and let~$\left(\Gamma_n\right)$ denote a sequence of subgroups
with the following properties:
\begin{itemize}
\item for each~$n$, $\Gamma_n\cong\mathbb Z^d$;
\item $\Gamma_1\subset\Gamma_2\subset\cdots$;
\item $\Gamma=\bigcup_{n\ge1}\Gamma_n$.
\end{itemize}
Let~$A_n$ denote the set of solutions to~\eqref{eq:Qcasethree}
with each~$\boldsymbol{q}_s^{(j)}\in\Gamma_n$
for which no subsum vanishes (thus~$A_n$
is a subset of the set
of values of~$j$ for which~\eqref{eq:Qcasethree}
holds, $A_1\subset A_2\subset\cdots$ and~$\bigcup_{n\ge1}A_n$ is the
set of all~$j$ for which~\eqref{eq:Qcasethree} holds).
We may assume without loss of generality that the
map
$$
j\mapsto\left({\boldsymbol{q}_2^{(j)}-
\boldsymbol{q}_1^{(j)}},\dots,
\boldsymbol{q}_r^{(j)}-
\boldsymbol{q}_1^{(j)}\right)
$$
is injective.
Since~$\Gamma_n$
is isomorphic to~$\mathbb Z^d$, Theorem~\ref{deepmagic}
applies to show that
\begin{equation}\label{eq:Qcaselevelnbound}
\left\vert A_n\right\vert\le\exp\left((6r)^{3r}(d+1)\right).
\end{equation}
Any finite bound in~\eqref{eq:Qcaselevelnbound} would
suffice to prove the main result in~\cite{MR1193598},
namely Theorem~\ref{Qcase} for~$\mathbb Z^d$-actions.
Here the additional uniformity of Theorem~\ref{deepmagic}
is needed:
The bound in~\eqref{eq:Qcaselevelnbound}
is independent of~$n$, so it follows that
equation~\eqref{eq:Qcasethree} holds without
vanishing subsum for only finitely many~$j$.
Thus there exists a set~$S\subsetneq\{2,\dots,r\}$
such that
\begin{equation}\label{eq:Qcasefour}
\sum_{s\in S}
(-a_s/a_1)\boldsymbol{u}^{\boldsymbol{q}_s^{(j)}-
\boldsymbol{q}_1^{(j)}}=0\mbox{ for infinitely many }j.
\end{equation}
The identity~\eqref{eq:Qcasefour} shows that~$\alpha$ is
not mixing on~$\vert S\vert<r$ sets.
Thus for any~$r<\infty$,
$$
\mathcal{M}(\alpha)\le r\implies\mathcal{M}(\alpha)<r,
$$
so
$$
\mathcal{M}(\alpha)>1\implies\mathcal{M}(\alpha)=\infty.
$$
This proves Theorem~\ref{Qcase} when~$I=P$ is prime.
Assume now that~$I$ is a radical ideal and that the system
corresponding to the module~$R_{\Gamma}/I$
is not mixing on~$r$ sets for some~$r>1$
but is mixing. As before,
this means
there is
a sequence of~$r$-tuples
$$
\left(
\boldsymbol{q}_1^{(j)},\dots,
\boldsymbol{q}_r^{(j)}
\right)_{j\ge1}\in\Gamma^r
$$
such
that~$\boldsymbol{q}_s^{(j)}-
\boldsymbol{q}_t^{(j)}\to\infty$ as~$j\to\infty$
for~$s\neq t$, for which
the equation
\begin{equation}\label{eq:Qcaseradical}
\boldsymbol{u}^{\boldsymbol{q}_1^{(j)}}a_1
+\dots+
\boldsymbol{u}^{\boldsymbol{q}_r^{(j)}}a_r=0\mbox{ for infinitely many }j
\end{equation}
holds in the ring~$R_{\Gamma}/I$.
Let
$$
U=\medspace\ll\negmedspace c,a_1,\boldsymbol{u}^{\boldsymbol{q}}-1\mid c\in
\mathbb Z\backslash\{0\},
\boldsymbol{q}\in\Gamma^d\negmedspace\gg
$$
where~$\ll\negmedspace A\negmedspace\gg$ denotes the multiplicative
group generated by~$A$.
This is a multiplicative set, and
we claim that~$U\cap I=\emptyset$.
If, for
some~$\boldsymbol{q}\in\Gamma^d$,
$\boldsymbol{u}^{\boldsymbol{q}}-1\in I$,
then~$\alpha_{R_{\Gamma}/I}$ is not mixing, which
is excluded by hypothesis.
Since~$I$ is radical, it follows
that~$(\boldsymbol{u}^{\boldsymbol{q}}-1)^m\notin I$ for all~$m\ge1$.
If
$$
(\boldsymbol{u}^{\boldsymbol{q}}-1)^mb\in I\mbox{ for
some }\boldsymbol{q}\in\Gamma^d,
m\in\mathbb Z\mbox{ and }b\notin I
$$
then we must
have~$(\boldsymbol{u}^{\boldsymbol{q}}-1)^mb^m\in I$. Since~$I$ is radical,
this implies
that~$(\boldsymbol{u}^{\boldsymbol{q}}-1)b\in I$, so in~$R_{\Gamma}/I$
$$
b+I=\boldsymbol{u}^{j\boldsymbol{q}}b+I\mbox{ for all }j\ge1,
$$
which again contradicts the assumption that~$\alpha_{R_{\Gamma}/I}$
is mixing.
By induction, no product of the form
$$
a_1^{\ell}
(\boldsymbol{u}^{\boldsymbol{q}_1}-1)^{m_1}\cdots
(\boldsymbol{u}^{\boldsymbol{q}_k}-1)^{m_k}
$$
can be in~$I$. Since~$R_{\Gamma}/I$ has
no additive torsion, it follows that~$U\cap I=\emptyset.$
By finding a maximal
ideal above~$I$ in the localization of~$R_{\Gamma}$ at~$U$,
we may find a prime ideal~$P\supset I$ with the
property that~$P\cap U=\emptyset$
(see~\cite[Prop.~2.11]{MR1322960}).
The equation~\eqref{eq:Qcaseradical}
drops via the map~$x\to x+P=\overline{x}$
to a non-trivial equation
\begin{equation}\label{eq:QcaseradicalinP}
\boldsymbol{u}^{\boldsymbol{q}_1^{(j)}}\overline{a_1}
+\dots+
\boldsymbol{u}^{\boldsymbol{q}_r^{(j)}}\overline{a_r}=0
\mbox{ for infinitely many }j,
\end{equation}
in which not all the coefficients have vanished.
It follows that the sequence
$$
\left(
\boldsymbol{q}_1^{(j)},\dots,\boldsymbol{q}_r^{(j)}
\right)_{j\ge1}\in\Gamma^r
$$
witnesses non-mixing on~$r$ sets in the prime system
corresponding to~$R_{\Gamma}/P$, and the argument above
shows that this is only possible
if~$\mathcal M(\alpha_{R_{\Gamma}/P})=1$.
By~\cite[Th.~1.6]{MR97c:28041}, this would require that there be a 
non-mixing element in~$\Gamma$, which is impossible
by the choice of~$U$.
\end{proof}

Theorem~\ref{Qcase} does not hold without the
assumption that the rational rank is finite -- see
Theorem~\ref{Qcrosscase}.
It also cannot hold without
the assumption of connectedness.
If the cyclic module~$R_{\Gamma}/I$ has 
additive torsion, then there is an element~$a+I$
and an integer~$k>0$ with~$ka\in I$.
For a sufficiently large~$n$,~$a\in R_{\Gamma_n}$,
and setting~$J=I\cap R_{\Gamma_n}$ induces an
inclusion
$$
R_{\Gamma_n}/J\subset R_{\Gamma}/I
$$
which
dualizes to show that the original action restricted
to~$\Gamma_n$ has a factor corresponding
to~$\mathbb Z[\Gamma_n]/J$ for some ideal~$J$
with~$J\cap\mathbb Z\neq\{0\}$. Since~$\Gamma_n\cong
\mathbb Z^d$, finite non-trivial order of mixing is
possible unless there are additional conditions on the
ideal by~\cite{MR97c:28041}.

One of the most striking features of~$\mathbb Z^d$-actions for~$d>1$
(as opposed to actions of~$\mathbb Z$) is that simple examples
may have order of mixing satisfying
$$
1<\mathcal M<\infty.
$$
This was
first pointed out by Ledrappier~\cite{MR80b:28030}, who
showed that the~$\mathbb Z^2$-action~$\alpha$
corresponding to the module~$\mathbb Z[u_1^{\pm1},u_2^{\pm1}]/
\langle2,1+u_1+u_2\rangle$
has
$$
\mathcal{M}(\alpha)=2;
$$
the paper~\cite{MR1971197} gives
a simple related
construction for any
specified order of mixing.
The same construction will give algebraic
$\mathbb Q^d$-actions for any~$d>1$ with
any specified order of mixing.

Theorem~\ref{Qcase}
shows that~$\mathbb Q^d$-actions on connected groups
also behave much like~$\mathbb Z^d$-actions.
Further evidence for the essential similarity of
algebraic~$\mathbb Q$-actions and~$\mathbb Z$-actions
is provided by the next result.

\begin{proposition}
Any mixing prime algebraic~$\mathbb Q$-action is mixing
of all orders.
\end{proposition}

\begin{proof}
Let the action correspond to the
module~$R_{\mathbb Q}/P$ for some prime ideal~$P$.
If the action is not mixing on~$r$ sets for some~$r\ge2$
then as before we find an equation
\begin{equation}\label{eq:Qcasetwooonedim}
{u}^{{q}_1^{(j)}}a_1
+\dots+
{u}^{{q}_r^{(j)}}a_r=0\mbox{ for infinitely many }j
\end{equation}
which holds in some field~$\mathbb K$ (not necessarily
of characteristic zero) with
$$
{q}_s^{(j)}-
{q}_t^{(j)}\to\infty\mbox{ in }\mathbb Q
\mbox{ as }j\to\infty
\mbox{ for }s\neq t.
$$
If~$P\cap\mathbb Z\neq\{0\}$ and~$P$
contains a non-trivial polynomial, then
the variable~$u$ satisfies an algebraic
equation over a finite field, so must be
a root of unity in~$\mathbb K$. This
precludes mixing.
So~$P=\langle p\rangle$ for some rational
prime~$p$, and the~$\mathbb Q$-action is
a full~$\mathbb Q$ shift on~$p$ symbols, so is mixing
of all orders.

If~$P\cap\mathbb Z=\{0\}$ then we are in
the setting of Theorem~\ref{Qcase}, which shows
that the system is mixing of all orders.
\end{proof}

\section{Actions of~$\mathbb Q_{>1}^{\times}$}

Actions of the much larger group~$\mathbb Q_{>1}^{\times}$
behave quite differently.
A situation that arises naturally for such actions is
one in which all the interesting
parts of the action are carried by a subgroup of~$\Gamma$
in the following sense. Suppose that the acting group is
a direct sum of two subgroups,
\begin{equation}\label{Gammasplits}
\Gamma=\Gamma'\oplus
\Gamma''
\end{equation}
and write~$\gamma=\gamma'+\gamma''$ for the
unique splitting of an element~$\gamma$ under~\eqref{Gammasplits}.
Suppose also that there is a~$\Gamma'$-invariant
splitting
\begin{equation}\label{splitspacemeasure}
(X,\mathcal B,\mu)\cong\prod_{\Gamma''}(X',\mathcal B',\mu'),
\end{equation}
and write elements of~$\prod_{\Gamma''}X'$ in the
form~$(x_{\beta})_{\beta\in\Gamma''}$.
Finally assume that
under these isomorphisms the transformation~$\alpha_{\gamma}$
on~$X$ is sent to the transformation
$$
(x_{\beta})_{\beta\in\Gamma''}\mapsto(y_{\beta})_{\beta\in\Gamma''}
$$
where
$$
y_{\beta}=\alpha_{\gamma'}(y_{\beta+\gamma''}).
$$
Write~$\overline{\alpha}$ for the image of the restriction of~$\alpha$
to~$\Gamma'$, acting on~$X'$.
That is, the subgroup~$\Gamma'$ acts by restricting
the original action, and the subgroup~$\Gamma''$ acts
as a full shift.
This kind of structure has been studied by Kaminski~\cite{MR629018}
in connection with the completely positive entropy property
for~$\mathbb Z^2$-actions. A \emph{cylinder set}
in~$\prod_{\Gamma''}(X',\mathcal B',\mu')$ is
a set of the form
$$
\{
(x_{\beta})_{\beta\in\Gamma''}\mid
x_{\beta_i}=a_i\mbox{ for }i\in F\}
$$
where~$F$ is a finite subset of~$\Gamma''$.

\begin{lemma}\label{independentsplitting}
$\mathcal M(\alpha)=\mathcal M(\overline{\alpha})$.
\end{lemma}

\begin{proof}
A sequence of elements in~$\Gamma'$ goes to infinity
if and only it goes to infinity in~$\Gamma$, so
$$
\mathcal M(\alpha)\le\mathcal M(\overline{\alpha}).
$$
Now let~$\left(\gamma^{(j)}_1,\dots,\gamma^{(j)}_r
\right)_{j\ge1}$ be a sequence of~$r$-tuples in~$\Gamma$ that
witnesses the statement~$\mathcal{M}(\alpha)<r$. This means
that
for every~$s\neq t$,
$$
\gamma^{(j)}_s-\gamma^{(j)}_t\rightarrow\infty\mbox{ as }
j\to\infty
$$
and -- without loss of generality --
there are cylinder sets~$A_1,\dots,A_r\in\mathcal B$ for which
$$
\mu\left(
\alpha_{\gamma^{(j)}_1}A_1\cap\dots\cap\alpha_{\gamma^{(j)}_r}A_r
\right)
\not\rightarrow\mu(A_1)\cdots\mu(A_r)\mbox{ as }
j\to\infty.
$$
If, for each~$s\neq t$,
$$
(\gamma^{(j)}_s)'-(\gamma^{(j)}_t)'\rightarrow\infty\mbox{ in }
\Gamma'\mbox{ as }
j\to\infty
$$
then the sequence~$\left((\gamma^{(j)}_1)',\dots,(\gamma^{(j)}_r)'
\right)_{j\ge1}$
witnesses that~$\mathcal M(\overline{\alpha})<r$,
so
$$
\mathcal M(\overline{\alpha})\le
\mathcal M(\alpha)
$$
as required.
If for some~$s$ and along some subsequence~$J$,
$$
(\gamma^{(j)}_s)''-(\gamma^{(j)}_t)''\rightarrow\infty\mbox{ in }
\Gamma''\mbox{ as }
j\to\infty\mbox{ in }J\mbox{ for all }t\neq s,
$$
then, since each~$A_i$ is a cylinder set,
$$
\mu\left(
\alpha_{\gamma^{(j)}_1}A_1\cap\dots\cap\alpha_{\gamma^{(j)}_r}A_r
\right)=
\mu(A_s)\times
\mu\left(
\bigcap_{t\neq s}\alpha_{\gamma^{(j)}_t}A_t
\right)\mbox{ for large }j\mbox{ in }J,
$$
so
$$
\mu\left(
\bigcap_{t\neq s}\alpha_{\gamma^{(j)}_t}A_t
\right)\not\rightarrow
\prod_{t\neq s}\mu(A_t)\mbox{ as }
j\to\infty\mbox{ in }J.
$$
Thus the sequence with term~$s$ omitted is
a witness to the statement that
$$
\mathcal M({\overline{\alpha}})<r-1,
$$
completing the proof.
\end{proof}

\begin{theorem}\label{Qcrosscase}
Let~$\alpha$ be an algebraic action of~$\mathbb Q^{\times}_{>0}$
corresponding to a module~$R_{\infty}/P$
with~$P\cap\mathbb Z=\{0\}$ and~$P$ a prime
ideal.
Then if~$P$ is finitely generated,
$$
\mathcal{M}(\alpha)>1\implies\mathcal{M}=\infty.
$$
If~$P$ is not finitely generated, the order
of mixing may be finite and positive.
\end{theorem}

\begin{proof}
View~$R_{\infty}$ as the ring~$\mathbb Z[u_1^{\pm1},\dots]$
of Laurent polynomials in countably many commuting variables
(this may be done via the map
$2\mapsto u_1$, $3\mapsto u_2$, $5\mapsto u_3$, and so on
through the primes).
If~$P$ is finitely generated, then
it has a set of generating polynomials which
only involve the variables~$u_1,\dots,u_d$ for some~$d\ge1$.
Then the action~$\alpha$ splits up as in Lemma~\ref{independentsplitting}
with~$\Gamma'\cong\mathbb Z^d$ and the action~$\overline{\alpha}$
being the~$\mathbb Z^d$-action corresponding to the
module~$\mathbb Z[u_1^{\pm1},\dots,u_d^{\pm1}]/P\cap
\mathbb Z[u_1^{\pm1},\dots,u_d^{\pm1}]$.
By Lemma~\ref{independentsplitting},
$$
\mathcal{M}(\alpha)=\mathcal{M}(\overline\alpha),
$$
and by~\cite{MR1193598},
$\mathcal{M}(\overline{\alpha})>1\implies
\mathcal{M}(\overline\alpha)=\infty$.

For the second statement, consider the natural
action of~$\mathbb Q^{\times}_{>0}$ on~$\widehat{\mathbb Q}$.
This is mixing, since the equation
$ax+b=0$ in~$\mathbb Q$ determines~$x$.
On the other hand, the action is not mixing on~$3$
sets since
$$
1\cdot(-1)+n\cdot(1)+(n-1)\cdot(1)=0
$$
for all~$n$ in~$\mathbb Q$,
and as elements of~$\mathbb Q^{\times}_{>0}$,
$n,n-1$ and~$\frac{n}{n-1}$ all go to infinity
and move apart as~$n\to\infty$.
\end{proof}

%\bibliographystyle{plain}
%\bibliography{references}

\end{document}